\documentstyle[amssymb,11pt]{article}

\date{March 14, 2014}
\title{Coding over Core Models\thanks{It is an honour for the authors
    to contribute this article in recognition of the 60th birthdays 
of Peter Koepke and Philip Welch.}}  
\author{Sy-David Friedman\thanks{The first author wishes to thank the
FWF for its support through Einzelprojekt P25671. He sees in Peter Koepke and Philip Welch fellow disciples of
our common mentor, Ronald Jensen.}\\
KGRC, Vienna
\and
Ralf Schindler\thanks{The result on pp.~\pageref{joint}f. was produced while the second author was visiting the
Erwin Schr\"odinger Institut, Vienna, in September 2013. He would like to 
thank
Sy Friedman and the other organizers of the ESI Set Theory Program for 
their
warm hospitality.}\\
WWU M\"unster
\and
David Schrittesser\thanks{The third  author wishes to thank Ralf Schindler for his support through SFB 878. He also wants to thank Sy Friedman 
and everyone at the KGRC for their
hospitality.}\\
KGRC, Vienna
}

\newtheorem{thm}{Theorem}

\newtheorem{cor}[thm]{Corollary}

\newcommand{\ZFC}{\hbox{ZFC}}
\newcommand{\ZF}{\hbox{ZF}}

\newcommand{\GCH}{\hbox{GCH}}

\newcommand{\ran}{\hbox{ran}}

\newcommand{\Ult}{\hbox{Ult}}
\newcommand{\bs}{\bigskip}
\newcommand{\noi}{\noindent}
\newcommand{\cof}{\hbox{cof }}

\begin{document}
\maketitle

Early in their careers, both Peter Koepke and Philip Welch made major
contributions to two important areas of set theory, core model theory (see
\cite{peter}) and
coding (see \cite{coding-book}), respectively. In this article we aim to survey some of the
work that has been done which combines these two themes, extending Jensen's
original Coding Theorem from $L$ to core models witnessing large
cardinal properties. 

\bs

The original result of Jensen can be stated as follows.

\begin{thm} (Jensen, see \cite{coding-book}) \label{jensen-coding} Suppose that $(V,A)$ is a
  transitive model of $\ZFC+\GCH$ (i.e., $V$ is a transitive model of $\ZFC+\GCH$ and
  replacement holds in $V$ for formulas mentioning $A$ as an
  additional unary predicate). Then there is a $(V,A)$-definable,
  cofinality-preserving class forcing $P$ such 
  that if $G$ is $P$-generic over $(V,A)$ we have:\\
(a) For some real $R$, $(V[G],A)\vDash \ZFC$ $+$ the universe is $L[R]$ and $A$ is definable with
parameter $R$.\\
(b) The typical large cardinals properties
consistent with $V=L$ are preserved from $V$ to $V[R]$: inaccessible, Mahlo, weak
compact, $\Pi^1_n$ indescribable, subtle, ineffable,
$\alpha$-Erd\H{o}s for countable $\alpha$.
\end{thm}

\begin{cor} \label{non-set-generic}
It is consistent to have a real $R$ such that $L$, $L[R]$ have the
same cofinalities but $R$ belongs to no set-generic extension of $L$.
\end{cor}

The theme of this article is to consider the following question: To
what extent is it possible to establish an analagous result when $L$
is replaced by a core model $K$ and the large cardinal properties in
(b) are strengthened to those consistent with $V=K$ (measurable,
hypermeasurable, strong, Woodin)?

\bs

A brief summary of the situation is as follows. Coding up to 
one measurable cardinal is unproblematic (see \cite{sy-code-meas}), although already in this
case there are some issues with condensation and the interesting new
phenomenon of ``ultrapower codings'' arises. At the level of hypermeasurable 
cardinals there are serious condensation issues which obstruct a fully
general result; nevertheless variants of Corollary
\ref{non-set-generic} can be established and
very special predicates $A$ as in Theorem \ref{jensen-coding} can be
coded (such as a generic for a Prikry product, see \cite{sy-mohammad}). In
addition, although one is able to lift enough of the total extenders on the
hierarchy of a core model witnessing hypermeasurability, it requires
extra effort to lift more than one total extender for the same critical point
(and it is not in general possible to lift all of the extenders (partial and total) 
on a fixed critical point $\kappa$ satisfying
$o(\kappa)=\kappa^{+++}$; we conjecture that this can be improved to
$o(\kappa)=\kappa^{++}$). At the level of Woodin cardinals, even 
Corollary \ref{non-set-generic} is not possible if the aim is to
lift all total extenders in a witness to Woodinness via the
``$A$-strong'' definition of this notion; however this obstacle is
removed by instead considering witnesses to the definition of
Woodinness in terms of ``$j(f)(\kappa)$ strength'' (see \cite{sy-glc}).

\bs

There are a number of applications of coding over core models. In
addition to those found in \cite{sy-book} based on Jensen's original
method, we mention two other examples.

\begin{thm} (Friedman-Schrittesser \cite{sy-david}) Relative to a Mahlo cardinal it is
  consistent that every set of reals in $L(\mathbb R)$ is Lebesgue measurable
  but some projective (indeed lightface $\Delta^1_3$) set of reals
  does not have the Baire property.
\end{thm}

\begin{thm} (Friedman-Golshani \cite{sy-mohammad}) Relative to a
  strong cardinal (indeed relative to a cardinal $\kappa$ that is
  $H(\kappa^{+++})$-strong) it is consistent to have transitive models
  $V\subseteq V[R]$ of $\ZFC$ where $R$ is a real, $\GCH$ holds in $V$
  and GCH fails at every infinite cardinal in $V[R]$. One can further
  require that $V$, $V[R]$ have the same cardinals.
\end{thm}

\noi
\emph{About Jensen coding}

\bs

To make what follows more intelligible it is worthwhile to first
review the case of Jensen coding. No matter how you look at it, even
this argument is complicated, although major simplifications can be
made if one assumes the nonexistence of $0^\#$ in the ground model
$V$. Our aim here however is not to delve into the fine points of the
proof (and in particular we will not reveal how the nonexistence of
$0^\#$ can be exploited), but rather to give the architecture of the
argument in order to facilitate a later discussion of generalisations.

\bs

For simplicity consider the special case where the cardinals of the
ground model $V$ are the same as those in $L$ and 
the ground model is $(L[A],A)$ where $A$ is a class of ordinals such
that $H(\alpha)=L_\alpha[A]$ for each infinite cardinal $\alpha$ (the
latter can be arranged using the fact that the GCH holds in $V$).

\bs

Coding is based on the method of almost disjoint forcing. Suppose that
$A$ is a
subset of $\omega_1$. Then we can code $A$ into  a real as follows: For
each countable ordinal $\xi$ attach a subset $b_\xi$ of $\omega$ (so that
the $b_\xi$'s are almost disjoint) and
force a real $R$ such that $R$ is almost disjoint from $b_\xi$ iff
$\xi$ belongs to $A$. Actually it is convenient to modify this to: 
$R$ almost contains $b_\xi$ iff $\xi$ belongs to $A$ (where
\emph{almost contains} means contains with only finitely many
exceptions). The conditions to achieve this are pairs
$(s,s^*)$ where $s$ is an $\omega$-Cohen condition (i.e. element of
$^{<\omega}2$) and $s^*$ is a finite subset of $A$; when extending to
$(t,t^*)$ we extend $s$ to $t$, enlarge $s^*$ to $t^*$ and insist that
if $s(n)$ is undefined but $t(n)$ equals $0$ then $n$
does not belong to $b_\xi$ for any $\xi$ in $s^*$. Then the generic
$G$ is determined by the union $G_0$ of the $s$ for $(s,s^*)$ in $G$ 
and we can take $R$ to be the set of $n$ such that $G_0(n)$ equals
$1$. The forcing has the ccc and ensures that $A$ belongs to $L[R]$
using the hypothesis $\omega_1=\omega_1^L$ to produce the $b_\xi$'s in
$L$ (and therefore also in $L[R]$).\footnote{In a more
general setting we have to worry about how to find the $b_\xi$'s in
$L[R]$. Jensen's trick to achieve this is to ``reshape'' $A$ into a
stronger predicate $A'$ with the property that any countable ordinal
$\xi$ is in fact countable in $L[A'\cap\xi]$; then after $R$ decodes
$A'\cap\xi$ it can find $b_\xi$ and continue the decoding. A
clever argument shows that such an $A'$ can be added over $L[A]$ by an
$\omega$-distributive forcing; when $A$ is not just a subset of
$\omega_1$ but a subset of some larger cardinal or even a proper class
of ordinals,
then the ``reshaping'' forcing must be woven into the coding forcing itself.}

\bs

There is nothing to stop us from coding a subset $A$ of $\omega_2$ into a
real in a similar fashion: First we use the hypothesis
$\omega_2=\omega_2^L$ to choose subsets $b_\xi$ of $\omega_1$ to set
up a forcing to code $A$ into a subset $B$ of $\omega_1$ via the
equivalence $\xi\in A$ iff $B$ almost contains $b_\xi$, and then we code
$B$ into a real as in the previous paragraph. It is pretty clear how to
do this for a subset of any $\omega_n$, $n$ finite. 

\bs

If we have a subset $A$ of $\aleph_\omega$ then we have to force subsets $A_n$ of
$\omega_n$ for each $n$ so that $A_n$ codes both $A_{n+1}$ and
$A\cap\omega_n$. At first this is confusing because there is no
``top'', i.e. no largest $n$ to begin with, but further reflection
reveals that there is no problem at all, as we don't need to know all of
$A_{n+1}$ to talk about conditions to add $A_n$. More precisely, a
condition $p$ will assign to each $n$ a pair $(s_n,s^*_n)$ so that
$s_n$ is an $\omega_n$-Cohen condition and $s^*_n$ is a size less than
$\omega_n$ subset of the set of $\xi$ such that $s_{n+1}(\xi)$ is defined
with value $1$. This makes sense even though $s_{n+1}$ is not defined
on all of $\omega_{n+1}$. We also insist that all of the $b_\xi$'s
consist of even ordinals and that each $A\cap\omega_n$ is coded into  
the union of the $s_n$'s using its values at odd ordinals. In the end $A$ gets coded into a
real and cofinalities are preserved since for any $n$ the forcing
factors into an $\omega_n$-closed forcing (the $n$-th upper part)
followed by an $\omega_n$-cc forcing (the $n$-th lower part).

\bs

Coding a subset $A$ of $\aleph_{\omega+1}$ into a real requires a new
idea. Actually by the previous paragraph it's enough to see how to
code $A$ into a subset of $\aleph_\omega$. Again we would like to assign a 
subset $b_\xi$ of $\aleph_\omega$ to each $\xi<\aleph_{\omega+1}$ and
then hope to force a subset $B$ of $\aleph_\omega$ which almost contains
$b_\xi$ iff $\xi$ belongs to $A$; how are we going to do
that? The conditions to add $B$ cannot be built from
``$\aleph_\omega$-Cohen conditions'' as this makes no sense for the
singular cardinal $\aleph_\omega$. Instead they should look like
conditions in the product of the $\omega_n$-Cohen forcings, i.e. of
the form $(s_n\mid n \in\omega)$ where each $s_n$ is an
$\omega_n$-Cohen condition (as in the previous paragraph but without
the ``restraints'' $s^*_n$). Actually it is very convenient to 
instead write $(s_{\omega_n}\mid n \in\{-1\}\cup\omega)$ where
$\omega_{-1}=0$ and $s_{\omega_n}$ is an $\omega_{n+1}$-Cohen
condition for each $n\geq -1$, and to think of
$s_{\omega_n}$ as an $\omega_{n+1}$-Cohen condition on the interval
$[\omega_n,\omega_{n+1})$ rather than on $\omega_{n+1}$, to separate
the domains of the different $s_{\omega_n}$'s. Thus 
the characteristic function of the generic subset $B$ of
$\aleph_\omega$ is the union of all of the $s_{\omega_n}$'s which appear in the
generic.  

\bs

As said above we'd like to
choose the $b_\xi$'s so that $\xi$
belongs to $A$ iff the generic subset $B$ of
$\aleph_\omega$ almost contains $b_\xi$ (i.e. contains $b_\xi$ with
a set of exceptions which is bounded in $\aleph_\omega$). This is
done using a \emph{scale}, i.e. a sequence $(f_\xi\mid \xi
<\aleph_{\omega+1})$ of functions in $\prod_{n\geq -1} [\omega_n,\omega_{n+1})$
which is cofinal mod finite. Then we take $b_\xi$ to be the range
of $f_\xi$. Again it is convenient to change notation: instead of
writing $f_\xi(n)$ we write $f_\xi(\omega_n)$. So the coding is: $\xi$
belongs to $A$ iff $G_{\omega_n}(f_\xi(\omega_n))=1$ for sufficiently large $n$, where
$G_{\omega_n}$ denotes the union of the $s_n$'s which appear in the
generic.

\bs

As we are using a scale we can arrange the following: if
$p=(s_{\omega_n}\mid -1\leq n<\omega)$ is a condition then for some ordinal
$|p|<\aleph_{\omega+1}$ called the \emph{height} of $p$, if $\xi$ is less than
$|p|$ then $\xi\in A$ iff
$s_{\omega_n}(f_\xi(\omega_n))=1$ for sufficiently large $n$ and if
$\xi$ is at least $|p|$ then 
$f_\xi(\omega_n)$ is not in the domain of $s_{\omega_n}$ for sufficiently large
$n$. In other words, $p$ already codes $A$ below $|p|$ but provides no
information about future coding on the interval $[|p|,\aleph_{\omega+1})$.
Notice the difference from the successor coding case: a single
condition will definitively code an initial segment of $A$, in the
sense that its values on a final segment of $b_\xi$ for $\xi$ in an initial segment of
$\aleph_{\omega+1}$ have already been fixed (restraints are not
needed). Of course no condition will code all of $A$, so this initial
segment of $A$ is proper. 

\bs

But how do we know that this coding of $A\subseteq \aleph_{\omega+1}$
into a subset of $\aleph_\omega$ preserves the cardinal
$\aleph_{\omega+1}$? For each $n$ we can factor the forcing as the
part $\geq\omega_n$ followed by the part below $\omega_n$, and as the
latter is a small forcing it causes no problems with
cardinal-preservation; so we want to show that the forcing
$\geq\omega_n$ (using conditions $p = (s_{\omega_k}\mid k\geq n)$) is
$\omega_{n+1}$-distributive, i.e. does not add new
$\omega_n$-sequences. For simplicity suppose that $n$ is $0$, so we
want to hit $\omega$-many open dense sets below any condition $p=
(s_{\omega_k}\mid k\geq 0)$. Here is the worry: maybe things are
going fine with the sequence $p=p_0\geq p_1\geq \cdots$ with
corresponding heights $|p_0|\leq |p_1|\leq\cdots$ so we can conclude
that the limit $p_\omega$ of the $p_n$'s will code $A$ up to the limit
$|p_\omega|$ of the $|p_n|$'s. But there is the danger that $p_\omega$ ``overspills''
in the sense that it already has assigned cofinally many values on
$b_\xi$ for some $\xi\geq |p_\omega|$. This 
unintended assignment may conflict with the desired coding of $A$ at the 
ordinal $\xi$.

\bs

The solution is to guide the construction using sufficiently elementary submodels
and to refine our concept of scale. Namely, when we build the $p_n$'s
we also build a definable $\omega$-chain of size $\aleph_\omega$ sufficiently elementary
submodels $M_0\prec M_1\prec\cdots$ of the universe which are
transitive below $\aleph_{\omega+1}$; we ensure that the $p_n$'s are
chosen from the $M_n$'s and have heights $|p_n|$ which interleave with
the ordinals $\gamma_n =M_n\cap\aleph_{\omega+1}$. The result is that the supremum of
the $|p_n|$'s is exactly $\gamma_\omega=M_\omega \cap\aleph_{\omega+1}$, where
$M_\omega$ is the union of the $M_n$'s. Now how does this help? The
point is that we can arrange for $p_\omega$, the limit of the $p_n$'s,
to be definable over
$M_\omega$ and therefore also over its transitive collapse $\overline
M_\omega$; if we can also arrange our scale so that $f_{\gamma_\omega}$ eventually
dominates any function in $\prod_n [\omega_n,\omega_{n+1})$ which is definable
over $\overline M_\omega$, then $p_\omega$ will leave a final segment
of the range of $f_{\gamma_\omega}$ untouched, as the sequence
$(|p_\omega(\omega_n)|\mid n\in\omega)$ is indeed definable over $\overline
M_\omega$. Finally, arranging our scale in this way is not a problem,
as $\overline M_\omega$ is an initial segment of $L$ which is so short that it still thinks that
$\gamma_\omega$ is a cardinal (it is the image of $\aleph_{\omega+1}$
under the transitive collapse of $M_\omega$) and we can define $f_\xi$
to eventually dominate any function in $\prod_n
[\omega_n,\omega_{n+1})$ which belongs to a model which still thinks
that $\xi$ is a cardinal ($f_\xi$ is defined using Skolem hulls inside
some big initial segment which sees that $\xi$ is not a cardinal).

\bs

The reason we discussed the fine point above about the coding of a
subset of $\aleph_{\omega+1}$ into $\aleph_\omega$ is to note that
there is some condensation involved (we needed that $\overline M_\omega$ is
an initial segment of our hierarchy). This is unproblematic for $L$
(and even for $L[U]$ where $U$ is a single normal measure) 
but is a serious problem for large core models. The
use of condensation is even more substantial when looking at
$\aleph_{\omega^2}$, where one needs to simultaneously consider transitive collapses
of unions of chains of sufficiently elementary submodels of any fixed
size $\aleph_{\omega\cdot n}$ and worry about their transitive
collapses being initial segments of the hierarchy. Indeed it is this
issue with condensation which obstructs a fully general coding result 
over core models as in Theorem \ref{jensen-coding}. Nearly all of the successes with
coding over core models are variants of the weaker Corollary \ref{non-set-generic}.

\bs

Now the fact that the strategy to code a subset of $\aleph_{\omega+1}$
into $\aleph_\omega$ fits so
nicely with the strategy to code a subset of $\aleph_\omega$ into a
real means that we can combine the two codings into a single coding of
a subset of $\aleph_{\omega+1}$ into a real. Thus a condition is a
function $p$ that for each finite $n$ assigns a pair
$(s_n,s^*_n)$ as in the latter coding so that in addition the sequence
of $s_n$'s is a condition in the former coding. For later use we
change notation slightly: the domain of $p$ consists of $0$
together with the $\omega_n$'s and for each $\alpha$ in the domain of
$p$, $p(\alpha)=(p_\alpha,p^*_\alpha)$ where $p_\alpha$ is an
$\alpha^+$-Cohen condition on the interval $[\alpha,\alpha^+)$
($0^+$ is taken to be $\omega$). And of course the restraint $p^*_\alpha$ is a size
at less than $\alpha^+$ subset of the set of $\xi$ such that
$p_{\alpha^+}(\xi)$ is defined with value $1$. We also require that
$p_\alpha$ codes $A\cap|p_\alpha|$ where the domain of $p_\alpha$ is
$[\alpha,|p_\alpha|)$, using its values at odd ordinals. Finally, for some
$|p|<\aleph_{\omega+1}$, if $\xi$ is less than $|p|$ then $\xi$
belongs to $A$ iff $p_\alpha(\eta)=1$  
for sufficiently large $\eta$ in $b_\xi=\ran(f_\xi)$ and when $\xi$ is
at least $|p|$, 
sufficiently large $\eta$ in $b_\xi$ lie outside
the domain of the $p_\alpha$'s.

\bs

This ends our introduction to Jensen coding. For arbitrary infinite
cardinals $\alpha$, the coding from a subset of $\alpha^{++}$ into a
subset of $\alpha^+$ is similar to the coding of a subset of
$\omega_1$ into a real and for arbitrary singular cardinals $\alpha$,
the coding of a subset of $\alpha^+$ into a subset of $\alpha$ is
similar to the above coding of a subset of $\aleph_{\omega+1}$ into a
subset of $\aleph_\omega$. The final case of the coding of a subset of
$\alpha^+$ into a subset of $\alpha$ for inaccessible $\alpha$ uses
either \emph{full support} and thereby resembles the singular coding, or
uses \emph{Easton support} and thereby resembles the successor
coding. In nearly all cases (including \cite{coding-book}) full
support is used (it faciliates the preservation of large
cardinals); Easton support coding is however needed in \cite{sy-david}.
The reason is that in \cite{sy-david}, we iterate Jensen coding to a length of $\kappa$,
at the same time collapsing everything below $\kappa$;
but we want to preserve $\kappa$ itself.
The usual strategy of ``reducing to the lower part'' fails below $\kappa$, since as we keep coding into $\omega$,
there are $\kappa$-many lower parts.
Instead, a much more complex argument is needed, in which there is no fixed height where we cut into ``upper'' and ``lower part'':
intuitively, we capture a given name by deciding it in different ways using larger and larger 
lower parts, catching our tail at an inaccessible below $\kappa$, where we will have looked at all the relevant
lower parts.
For this to work, supports must be bounded below inaccessibles 
(we also have to assume $\kappa$ is Mahlo).

\bs

\noi
\emph{One measurable cardinal}

\bs

Suppose that there is a measurable cardinal $\kappa$ in $V$. Can we
code $V$ into a real $R$ preserving the measurability of $\kappa$? 

\bs

Of course the model that results after coding into $R$ cannot be
$L[R]$, but it could be $L[U^R,R]$ where $U^R$ is a normal measure on
$\kappa$ extending a given normal measure $U$ on $\kappa$ in
$V$. As alluded to above there are serious issues with condensation when
coding over core models and
for this reason we'll only discuss here how to establish a version of
Corollary \ref{non-set-generic}: It is possible to force a real $R$
over $L[U]$ which preserves cofinalities, is not set-generic over
$L[U]$ and preserves the measurability of $\kappa$. Even in this
special situation it is very helpful (and essential for 
further generalisations) to use a hierarchy for $L[U]$ with good
condensation properties, which we write as $L[E]$. Note that the
$L[U]$-hierarchy does not obey even the weakest of consequences of
condensation, the property that subsets of an infinite cardinal
$\alpha$ appear in the hierarchy at a stage before $\alpha^+$. The
$L[E]$ hierarchy inserts ``partial measures'' which ensure this
property and more without altering the model: $L[E]=L[U]$. The measure
$U$ (or something very close to it) is placed on the $L[E]$ hierarchy
at an appropriate stage between $\kappa^+$ and $\kappa^{++}$, its
\emph{index} on the $L[E]$-hierarchy, and there will be many
approximations to it placed on the hierarchy at indices cofinal in any
uncountable cardinal up to and including $\kappa^+$. 

\bs

So proceed now to form conditions $p$ in $L[E]$ which resemble the
coding conditions from Jensen coding: For $\alpha$ either $0$ or an
infinite cardinal, $p(\alpha)$ is a pair $(p_\alpha,p^*_\alpha)$ where
$p_\alpha$ is an $\alpha^+$-Cohen condition on $[\alpha,\alpha^+)$ and $p^*_\alpha$ is a size
at most $\alpha$ set of $\xi$ such that $p_{\alpha^+}(\xi)=1$. Also
for limit cardinals $\lambda$ we have a scale
$(f_\xi\mid\xi\in[\lambda,\lambda^+))$ of functions in
$\prod_{\alpha^+<\lambda} [\alpha^+,\alpha^{++})$ and for
$\xi<|p_\lambda|$, $p_{\lambda}(\xi)= 1$ iff
$p_{\alpha^+}(f_\xi(\alpha^+)) =1$ for sufficiently large
$\alpha^+<\lambda$. And $p\upharpoonright\lambda$ does not interfere
with future coding on $[|p_\lambda|,\lambda^+)$ in the sense that for
$\xi\ge |p_\lambda|$, $p_{\alpha^+}(f_\xi(\alpha^+))$ is not defined
for sufficiently large $\alpha^+<\lambda$. The previous applies both to
inaccessible and singular limit cardinals $\lambda$. 

\bs

Now we need a strategy for showing that this forcing preserves the
measurability of $\kappa$. It is best to think of measurability in
terms of embeddings: In the ground model $V=L[U]=L[E]$ there is an
elementary embedding $j:V\to M=\Ult_U$ with critical point $\kappa$, derived
from the ultrapower given by $U$. The hierarchy provided by $E$ is 
defined so that we have $j:L[E]\to L[E^*]$ where $E$, $E^*$ agree up to
the index of $U$ (an ordinal between $\kappa^+$ and $\kappa^{++}$);
for the present discussion we only need to know that this agreement
persists at least up to the $\kappa^{++}$ of $M$, the
ultrapower of $V$ by $U$. This has the important consequence that our
coding forcing $P$ agrees with $P^*=j(P)$, the coding forcing of the
ultrapower $M$, up to the $\kappa^{++}$ of $M$. More precisely, a
function $p$ defined at $0$ together with the infinite cardinals
$\le\kappa^+$ such that $p(\alpha)=(p_\alpha,p^*_\alpha)$ for each
$\alpha$ and $p^*_{\kappa^+}=\emptyset$ 
belongs to $P^*$ iff it belongs to $P$, $|p_{\kappa^+}|$
is less $(\kappa^{++})^M$ and $p^*_\kappa$ is a subset of 
$(\kappa^{++})^M$.\footnote{This may not be entirely clear, as $V$
  has more subsets of $\kappa^+$ than $M$. However the coding is
  defined so that $p_{\kappa^+}$ will belong to $M$ provided its
  length is less than the $\kappa^{++}$ of $M$.} 

\bs

Now Silver taught us that if we want to preserve the measurability of
$\kappa$ we should lift the embedding $j:V\to M$ to an embedding
$j^*:V[G]\to M[G^*]$ where $G^*$ is generic over $M$ for $P^*=j(P)$.
The key is to choose $G^*$ to contain the pointwise image $j[G]$ of
$G$ as a subclass. There are many examples of such liftings in the
context of reverse Easton forcing, where there are typically many
choices for $G^*$. But notice that with coding there is only one
candidate for $G^*$, the $P^*$-generic coded into the same real $R$ that
codes $G$. This is because $j^*(R)$ will equal $R$ for any possible
lifting $j^*$ of $j$ to $V[G]$. 

\bs

Of course our desired generic $G^*$ must include the image $j(p)$ of
any condition $p$ in $G$; it would be ideal if $G^*$ were 
simply generated by these conditions in the sense that $G^*$ is
obtained as the class of all conditions extended by a condition in
$j[G]$. This will however not be the case and it is instructive to see
why not. 

\bs

For $G^*$ to be generic it must intersect all $L[E^*]$-definable dense classes
$D$ on the forcing $P^*$. 
As $L[E^*]$ is the ultrapower of $L[E]$ by
the measure $U$ we can write $D$ as $j(f)(\kappa)$ for some definable function
$f$ with domain $\kappa$ in $L[E]$ so that $f(\alpha)$ is dense on
$P$ for each $\alpha$. Now our coding forcing $P$ satisfies the
following useful form of ``diagonal distributivity'': We say that
a subclass $D$ of $P$ is \emph{$\gamma$-dense} for a cardinal $\gamma$
if any condition in $P$ can be extended into $D$ without changing its
values below $\gamma$.  Now suppose that
$f(\alpha)$ is $\alpha^+$-dense for each cardinal $\alpha<\kappa$ and $p$ is a
condition. Then $p$ has an extension $q$ which meets 
(i.e.~extends an element of) each
$f(\alpha)$. It follows that some condition $p$ in
$G$ meets each $f(\alpha)$ and therefore on the
ultrapower side, $j(p)$ will meet $j(f)(\kappa)=D$ provided $D$ is
$\kappa^+$-dense on $P^*$. In particular 
this means that the $j(p)$ for $p$ in $G$ will
indeed provide us with a generic for the forcing $P^*$ above
$\kappa^+$, i.e. a generic subset of the $\kappa^{++}$ of $L[E]$ that in
turn codes an entire generic class for the forcing $P^*$ above
$\kappa^+$. As the embedding $j$ is the identity below $\kappa$,
$j[G]$ also provides us with a generic below $\kappa$ and indeed a
generic subset $G_\kappa$ of $\kappa^+$, as this is coded in both $L[E]$ and $L[E^*]$ 
into the generic below $\kappa$ in the same way.

\bs

So $j[G]$ in fact gives us a 
subset $G^*_{\kappa^+}$ of $(\kappa^{++})^M$ which codes an entire $P^*$-generic 
above $(\kappa^{++})^M$, as well as a 
subset $G_\kappa=G^*_\kappa$ of $\kappa^+$ which 
is generically coded (in both the $P$ and $P^*$ forcings) into a real;
what is missing is to
ensure that $G_\kappa$, which generically codes $G_{\kappa^+}$ over
$V[G_{\kappa^+}]$, also 
generically codes $G^*_{\kappa^+}$ over $M[G^*_{\kappa^+}]$. We have
to fit the ``ultrapower coding'' of $G^*_{\kappa^+}$ into $G_\kappa$ together 
with the ``$V$-coding'' of $G_{\kappa^+}$ into $G_\kappa$, in order to
produce the desired $P^*$-generic $G^*$.

\bs

It is tempting now to make use of the fact that $V=L[E]$ and
$M=L[E^*]$ actually agree up to $(\kappa^{++})^M$ in the sense that
the hierarchies given by $E$ and $E^*$ are the same up to that point. 
Indeed it is natural to expect that $G_\kappa$ will
generically code $G^*_{\kappa^+}$ using $E\upharpoonright 
(\kappa^{++})^M$, since it generically codes $G_{\kappa^+}$ using $E$
and $E^*\upharpoonright (\kappa^{++})^M$ is an initial segment of $E$.
This is encouraging, however it leads to a contradiction,
as what $G_\kappa$ codes below the ordinal $(\kappa^{++})^M$ using $E$
is $G_{\kappa^+}$ restricted to this ordinal, an element of $V$, whereas
what we want $G_\kappa$ to code over $M$, namely
$G^*_{\kappa^+}$, cannot be an element of $V$ (else both 
$G^*_{j(\kappa)}$ and its preimage $G_\kappa$ would belong to $V$, 
reducing our class-forcing to a set-forcing). 

\bs

Thus we need a different approach, in which the codings over $L[E]$
and $L[E^*]$ do not agree at $\kappa^+$, in the sense that the generic
subset $G_\kappa$ of $\kappa^+$ codes the generic subset
$G_{\kappa^+}$ of $\kappa^{++}$ using $E$ in a way which accomodates, but
differs from, the way it codes the generic subset
$G^*_{\kappa^+}$ using $E^*$. The solution is this: When defining conditions
$p(\kappa)=(p_\kappa,p^*_\kappa)$ to almost disjoint code
$p_{\kappa^+}:[\kappa^+,|p_{\kappa^+}|)\to 2$ we use sets $b_\xi$ for
$\xi<\kappa^{++}$ as before to ensure that $p_{\kappa^+}(\xi)=1$ iff
$p_{\kappa}(\delta)=1$ for sufficiently large $\delta\in b_\xi$;
however we additionally have sets $b^*_\xi$ for $\xi <
(\kappa^{++})^M$ to ensure that for $\xi<(\kappa^{++})^M$, $j(p)_{\kappa^+}(\xi)=1$ iff
$p_{\kappa}(\delta)=1$ for sufficiently large $\delta\in
b^*_\xi$. Thus there are two codings taking place simultaneously, one
of $p_{\kappa^+}$ and the other of $j(p)_{\kappa^+}$, with two
different forms of restraint. To avoid
conflicts between these codings we choose the $b_\xi$'s to be very
``thin'' making use of the measure $U$. We choose a scale
$(f_\xi\mid\xi\in[\kappa^+,\kappa^{++}))$ of functions from $\kappa^+$
to $\kappa^+$ so that the least function $f_{\kappa^+}$ of this scale 
eventually dominates all functions from $\kappa^+$ to $\kappa^+$ in
$M=L[E^*]$; this is possible as there are only $\kappa^+$-many such
functions in $M$. The net effect is that the resulting subset $G_\kappa$
of $\kappa^+$ which is generic over $L[E]$ will also be generic over
$L[E^*]$, as the thinness of the sets $b_\xi$ allows us to show that
conditions can be extended to meet the necessary dense sets from
the $L[E^*]$ coding without conflicting with the restraint imposed by
the $b_\xi$ for $\xi$ in $p^*_{\kappa^+}$. 

\bs

\noi
\emph{Measures of higher order}

\bs

Suppose now that we are in a ``Mitchell model'' $L[E]$ where we now
have two normal measures $U_0,U_1$ on $\kappa$ with $U_0$ below $U_1$ in the
Mitchell order. Thus $U_0$ belongs to the ultrapower of $V$ by the
measure $U_1$. Can we create a real which is class-generic but not
set-generic lifting both of the measures $U_0$ and $U_1$?

\bs

It is convenient to reformulate the situation of the last section (with
a single measure $U$) as follows. Recall that at $\kappa^+$ we have two
codings, that of $j(p)_{\kappa^+}$ into $G_\kappa\subseteq\kappa^+$
over the ultrapower $\Ult_U$ of $V$ by $U$, and the other of
$p_{\kappa^+}$ into $G_\kappa$ over $V$. As the latter coding takes
place ``above'' the former (ultrapower) coding, it is natural to think
of the $\kappa^{++}$-Cohen condition
$p_{\kappa^+}:[\kappa^+,|p_{\kappa^+}|)\to 2$ in two parts: there is
$p_{\kappa^+}$ on $[\kappa^+,|j(p)_{\kappa^+}|)$ coinciding with
$j(p)_{\kappa^+}$ and then $p_{\kappa^+}$ on
$[(\kappa^{++})^{\Ult_U},|p_{\kappa^+}|)$, which is coded using the $b_\xi$'s
which ``lie above'' the ultrapower $\Ult_U$. In this way there is in a
sense just one coding, which uses restraints from $\Ult_U$ below the
$\kappa^{++}$ of $\Ult_U$ and restraints from $V$ between the
$\kappa^{++}$ of $\Ult_U$ and the real $\kappa^{++}$. In fact
$p\upharpoonright\kappa$ is responsible for the coding below
$\kappa^{++}$ of $\Ult_U$ (via
the embedding $j$) and $p_{\kappa^+}$  is responsible for the coding
above. But notice that viewed this way, the domain of $p_{\kappa^+}$
is no longer an interval, but the union of two intervals, namely 
$[\kappa^+,|j(p)_{\kappa^+}|)$ and
$[(\kappa^{++})^{\Ult_U},|p_{\kappa^+}|)$. So $p_{\kappa^+}$ is what
one might call a ``perforated string''.

\bs

Now let's return to the more complex case of two measures $U_0$,
$U_1$. At $\kappa^+$ the strings are doubly-perforated, as their
domains consist of the union of three intervals:
$[\kappa^+,|j_{U_0}(p)_{\kappa^+}|)$,
$[(\kappa^{++})^{\Ult_{U_0}},|j_{U_1}(p)_{\kappa^+}|)$ and
$[(\kappa^{++})^{\Ult_{U_1}},|p_{\kappa^+}|)$. For 
cardinals $\bar\kappa$ of Mitchell order $0$ (i.e. carrying only
normal measures concentrating on non-measurables), strings at
$\bar\kappa^+$ will only be singly-perforated and at non-measurables
we return to non-perforated strings. The situation is similar, but
more complicated, when dealing with measurable $\kappa$ of Mitchell
order less than $\kappa^{++}$ (the ``real'' coding takes place above
the supremum of the $(\kappa^{++})^{\Ult_U}$ for $U$ a normal measure
on $\kappa$ on the $L[E]$ hierarchy). 

\bs

But if we go as far as $o(\kappa) = \kappa^{++}$, where
$(\kappa^{++})^U$ can be arbitrarily large in $\kappa^{++}$ for 
measures $U$ on $\kappa$, and wish to lift all of these measures,
then we have a problem, as it seems that 
there is no longer room to code, as the entire interval
$[\kappa^+,\kappa^{++})$ has been covered with ultrapower codings
which must be respected. 
In fact, we cannot expect to add a class-generic real which is not
set-generic but lifts all extenders, partial and total, 
when $o(\kappa)=\kappa^{+++}$:

\bs
\label{joint}
To see this, work in  $K=L[E]$ and assume that $o(\kappa) = \kappa^{+++}$,
$\kappa$ is the largest measurable cardinal, but $L[E]$ is also closed under sharps.
Let $S^0$, $S^1$ be a canonical partition of $S=\{ \xi < \kappa^{+++} \;|\; \cof(\xi) = \kappa^{++}\}$ into stationary sets.
Fix $i\in\{0,1\}$.
Let $T^i = S^i \cup \{ \xi < \kappa^{+++} \;|\; \cof(\xi) < \kappa^{++}\}$.
For later use, let $T^i_\lambda$ denote the set defined just like
$T^i$ 
but with $\kappa$ replaced by $\lambda$, for regular $\lambda$.
Assume now that there is a forcing $P^i$ in $K$ 
such that in the $P^i$-generic extension there exists
a real $r^i$ with the following property:

\bs

\noi
$(*)$ For all $\xi$ such that 
$K|\xi \vDash$``$\ZF^-$, $\lambda$ is the largest measurable cardinal and  $o(\lambda)=\lambda^{+++}$'', there is a club through $(T^i_\lambda)^{K|\xi}$
in $K|\xi[r^i]$.
Moreover every (partial or total) extender on $E$ lifts to $K[r^i]$.

\bs

By the last sentence we mean that every iteration tree on $K$ can be lifted to one on $K[r]$.
Now force with $P = P^0 \times P^1$, 
obtaining reals $r^0, r^1$ as above
(note that $P$ collapses $\kappa^{+++}$, but this is irrelevant), 
and let $g$ be generic for the collapse of 
$\kappa^{++}$ to $\omega$ over $K[r^0, r^1]$.
By an unpublished construction of Woodin  (see  \cite{doebler}), $K[g]$ is $\Sigma^1_4$-correct
in $K[g, r^i]$; this makes vital use of the fact that enough extenders from $K$ 
lift to $K[r^i]$ (cofinally many total extenders suffice here).
But the following $\Pi^1_3$ statement $\Psi^i(r^i)$ holds in $K[g,r^i]$, where we let $\Gamma(\lambda)$ denote the theory ``$\ZF^- \wedge\lambda$ is the largest measurable cardinal and  $o(\lambda)=\lambda^{+++}$'':
\bs

\noi
$(**)$ Every countable mouse $M_0$ such that $M_0\vDash \Gamma(\lambda)$ 
has a simple countable 
iterate $M_1$ such that $r^i$ lifts all extenders on the $M_1$-sequence,
and if $M_0$ is a countable mouse such that $M_0\vDash \Gamma(\lambda)$ 
and $r^i$ lifts all extenders on the $M_0$-sequence, then 
there is a club through $(T^i_{\lambda})^{M_0}$ in $M_0[r^i]$.

\bs
That  this statement is $\Pi^1_3$ boils down to the fact that being a mouse, in our setting, is $\Pi^1_2$.
By a simple iterate we mean that there are no drops.
To see that the first part of $(**)$ 
holds of $r^i$, let $M_0$ be a countable mouse.
Co-iterate $M_0$ with $K$ until you reach $M_1 \triangleleft K'$, where $K'$ is an iterate of $K$.
Since $r^i$ lifts all extenders on $E$, we can push forward $(*)$ to $K'$ in the sense that $(*)$ holds with initial segments of $K$ replaced by those of 
$K'$,  
so $M_1$ is as desired. If $M_1$ is not countable, then by taking a countable hull
we will obtain a {\em countable} iterate of $M_0$ which is as desired.
To see the second part of $(**)$, we may argue in a similar fashion, this
time co-iterating $M_0[r^i]$ with $K[r^i]$.

\bs

By $\Sigma^1_4$-correctness, we find $s^0, s^1$ in $K[g]$ such that 
$K[g]\vDash \Psi^0(s^0)\wedge\Psi^1(s^1)$.

\bs

We claim that for $i \in \{ 0 , 1 \}$,
$K$ has a simple iterate $K'$ such that $s^i$ lifts all extenders on the 
$K'$-sequence and there is a club through $(T_{\kappa'}^i)^{K'}$, where 
$\kappa'$ is the image of $\kappa$ under the iteration map.
Otherwise this is false for some $K|\theta \models \Gamma(\lambda)$, 
and we may pick some
cardinal
$\Omega >\theta$ and some 
$\sigma \colon {\bar K}[s^i] \rightarrow K|\Omega[s^i]$ with $\theta \in {\rm ran}(\pi)$
and ${\bar K}[s^i] \in K[g]$ is transitive and countable in $K[g]$.
Let $h \in K[g]$
be ${\rm Col}(\omega,\theta)$-generic over ${\bar K}[s^i]$. By our 
hypothesis that $L[E]$ is closed under sharps, 
${\bar K}[s^i][h]$ will be $\Sigma^1_2$-correct in $K[g]$.
This means that if we look at the family ${\cal F}$ of all $M \in {\bar K}[s^i][h]$
which in ${\bar K}[s^i][h]$ are countable iterates of ${\bar K}|\theta$,
then using $(**)$ densely many $M \in {\cal F}$ will be such that $r^i$ 
lifts all extenders on the $M$-sequence, and there is a club through
$(T_\lambda^i)^M$ in $M[r^i]$ (where $M \models \Gamma(\lambda)$).
But then if ${\tilde M}$ is the direct limit of all mice in ${\cal F}$,
then ${\tilde M} \in {\bar K}[s^i]$ by the homogeneity of  ${\rm Col}(\omega,\theta)$,
$r^i$ 
lifts all extenders on the ${\tilde M}$-sequence, and there is a club through
$(T_\lambda^i)^{\tilde M}$ in ${\tilde M}[r^i]$ (where ${\tilde M} \models \Gamma(\lambda)$). Moreover, in ${\bar K}[s^i]$, ${\tilde M}$ can be absorbed by
an {\em iterate} of ${\bar K}|\theta$ which has the same properties.
But the elementarity of $\sigma$ then yields a contradiction.

\bs

We first find a simple iterate $K'$ of $K$ as in $\Psi^0(s^0)$.
There is a club $D^0$ through $j_0(T^0)$ in $K'[s^0]$, where
$j_0\colon K \rightarrow K'$ is the resulting embedding. 
We can also find a further simple iterate $K''$ of $K'$ which is a witness to $\Psi^1(s^1)$.
We have $j_1\colon K' \rightarrow K''$ and 
since $s^0$ lifts all extenders on $E^{K'}$, we can assume
$j_1\colon K'[s^0] \rightarrow K''[s^0]$.
Thus there is a club $D^1$ through $j_0(T^1)$ in $K''[s^1]$, where $j = j_1\circ j_0$.
Let $\bar D^0 = j_1(D^0)$, noting this is a club through $j(T^0)$, and
let $D$ denote the limit points of $j^{-1}[\bar D^0\cap D^1]$.
Obviously, $D$ is club in $\kappa^{+++}$ and $D\in K[g]$.
As $S$ remains stationary in $K[g]$, we can find $\eta \in S\cap 
D$.
Since $j$ is continuous at points whose $K$-cofinality is greater than $\kappa$, we have $j(\eta) \in \bar D^0\cap D^1$, and by elementarity $j(\eta) \in j(S)$.
This is a contradiction, since $j(S) \cap j(T^0) \cap j(T^1) =\emptyset$, finishing the argument that no real $r^i$ as in $(*)$ can exist.

\bs

We now argue that the problem with $(*)$ lies not in its first sentence,
the coding part, but in its second sentence, the lifting of extenders. 
Consider the following weakening of $(*)$:

\bs

\noi
For all $\xi$, \emph{if}
$K|\xi \vDash$ ``$\lambda$ is the greatest measurable and  $o(\lambda)=\lambda^{+++}$''
and $K|\xi[r^i] \vDash$ ``$\ZF^-$ and $\lambda^{+++}=(\lambda^{+++})^{K|\xi}$'', \emph{then} there is a club through $(T^i_\lambda)^{K|\xi}$
in $K|\xi[r^i]$,

\bs
We can produce a real $r^i$ satisfying this by first shooting a club through $T^i$ and then forcing to code it with localization (using a core model analogue of David's trick; see \cite[theorem 6.18]{sy-book}). 
The club is added by a $\kappa^{+++}$ distributive forcing (a forcing adding no new $\kappa^{++}$-sequences) of $K$;
then, the condensation provided by $K$ suffices for the distributivity of the second forcing, as we can take Skolem hulls in $K$ (as in \cite{sy-mohammad}).
In fact, if we weaken $(*)$ by just dropping the last requirement (i.e. that all extenders lift), we obtain a statement which should be forceable using a core model analogue of \emph{strong coding} (see \cite{sy-sc}). 
For these reasons we believe that the problem with $(*)$ lies with its second assertion, that all extenders lift.

\bs

In fact we conjecture that it is not possible to add a class-generic
real which is not set-generic while lifting all 
normal measures on a measurable $\kappa$ of order $\kappa^{++}$; 
in fact we conjecture that it is not even possible to do this while lifting 
all normal measures in a ``cofinal'' collection ${\cal S}\in L[E]$ of
such measures. (By ``cofinal'' we mean that the ordinals
$(\kappa^{++})^{\Ult_U}$ for $U$ in $\cal S$ are 
cofinal in $\kappa^{++}$.)
But this does \emph{not} mean that we cannot preserve the property $o(\kappa)=
\kappa^{++}$! As we'll see in the next section, it is possible to
preserve cases of hypermeasurability, which in turn implies that 
the set of normal measures that are lifted is cofinal; it is not
however clear that this set can contain a cofinal subset in the ground
model.

\bs

\noi
\emph{Hypermeasurables}

\bs

Can we preserve stronger forms of measurability? Suppose that $\kappa$
is hypermeasurable in $V=L[E]$ in the sense that some total extender $F$ on
the $E$-sequence with critical point $\kappa$ witnesses that $\kappa$
is $H(\kappa^{++})$-strong, i.e. the ultrapower $j_F:V\to M$ has the
property that $H(\kappa^{++})$ is contained in $M$. (This is the same
as saying that $F$ is indexed past $\kappa^{++}$ in the
$E$-hierarchy.) Can we add a real which is class-generic but not
set-generic and lifts $F$? 

\bs

Again we want to set up our conditions so that the embedding $j=j_F$
can be lifted to $V[G]$. This time we have that the union of the
$j(p)_{\kappa^{++}}$ is not in $V$ yet like $p_{\kappa^{++}}$ must be
coded into the same subset $G_{\kappa^+}$ of $\kappa^{++}$. As in the
one measure case this can be resolved by starting the coding of
$p_{\kappa^{++}}$ above
$(\kappa^{+++})^{\Ult_F}$, below which the former coding takes
place. But we have a new problem: the union of the $j(p)_{\kappa^+}$
would appear to not belong to $V$ and as $V$ and $\Ult_F$ completely agree below
$\kappa^{++}$, the set $G_{\kappa^+}$ will code it in exactly the same
way as it codes $G_{\kappa^{++}}$. This is a serious obstacle and the
only way around it is to thin out the coding conditions to ensure
that in fact the $j(p)_{\kappa^+}$ will be empty for each condition
$p$.

\bs

To ensure the latter we require that for any condition $p$ 
there is a closed unbounded subset $C$ of $\kappa$ such that for
inaccessible $\alpha$ in $C$, $p_{\alpha^+}$ is the empty string. This
ensures that $j(p)_{\kappa^+}$ will also be the empty string. The
price one pays for this is that we only have a weaker form of diagonal
distributivity: If $f(\alpha)$ is $\alpha^{++}$-dense for each
cardinal $\alpha<\kappa$ then any condition can be extended to 
meet each $f(\alpha)$. This only ensures that the
pointwise image $j[G]$ will generate a generic over the ultrapower
$\Ult_F$ above $\kappa^{++}$, coded into the subset $G^*_{\kappa^{++}}$
of $(\kappa^{+++})^M$ consisting of the union of the
$j(p)_{\kappa^{++}}$ for $p$ in $G$. $G$ provides a generic below $\kappa^{++}$
and now the task is to ensure that $G_{\kappa^+}$ will code not only
$G_{\kappa^{++}}$ but also $G^*_{\kappa^{++}}$. This is dealt with as
in the one measure case, by starting the former coding ``above'' the
latter, making use of an appropriate scale. 

\bs

For a stronger total extender (of sucessor cardinal strength) the pattern is similar: Thin out the
conditions to guarantee that $j(p)_\alpha$ is empty for cardinals
$\alpha$ strictly between $\kappa$ and the strength of the
total extender. At the strength there are two codings which must be
performed simultaneously, one over $V$ and the other over the
ultrapower. Conflicts between these codings are avoided by allowing the
$V$-coding to make use of the total extender $F$ when defining the coding
sets $b_\xi$ via an appropriate scale.

\bs

The above ideas are sufficient to lift a class $\cal S$ of total extenders
(each of successor cardinal strength) which
is \emph{bounded} (the set of $(\alpha^+)^{\Ult_F}$
for total extenders $F$ in $\cal S$ of strength exactly $\alpha$ is bounded in $\alpha^+$
for each cardinal $\alpha$) and \emph{uniform (or coherent)} (if $F$ belongs
to $\cal S$ then $j_F({\cal S})$ agrees with $\cal S$ below the index of $F$ in the
$L[E]$-hierarchy), provided that in $L[E]$ no inaccessible $\alpha$ is the
stationary limit of cardinals which are strong up to $\alpha$. This
yields a version of Corollary \ref{non-set-generic} up to the level of
a proper class of strong cardinals, but handling a stationary-limit of
strong cardinals requires new ideas.

\bs

\noi
\emph{Woodin cardinals}

\bs

As coding makes heavy use of condensation it is only reasonable to
consider ground models for which a suitable core model theory is available,
currently up to the level of Woodin cardinals. 

\bs

Recall that $\delta$ is \emph{Woodin} if for each $A\subseteq\delta$
there is a $\kappa<\delta$ which is \emph{$A$-strong in $\delta$},
i.e. the critical point of embeddings $j:V\to M$ such that $j(A)$
agrees with $A$ up to $\gamma$, for each $\gamma<\delta$. At first it
appears that this indicates the end of the coding method, as Woodin
proved the following (see \cite{sy-glc} and \cite[theorem 7.14]{hbst:steel}): If $\cal S$ is a set of total extenders in $V$ sufficient
to witness Woodinness in this sense and $R$ is a real such that each
total extender in $\cal S$ lifts to $V[R]$, then in fact $R$ is generic over $V$
for a ($\delta$-cc) forcing of size $\delta$. So there appears to be no version
of Corollary \ref{non-set-generic} in the context of a Woodin
cardinal.

\bs

But actually there is another definition of Woodin cardinal with a
different notion of witness: $\delta$ is Woodin if for each
$f:\delta\to\delta$ there is a $\kappa<\delta$ closed under $f$ which
is \emph{$f$-strong}, meaning that some embedding $j:V\to M$ with
critical point $\kappa$ is $j(f)(\kappa)$-strong (i.e.,
$H(j(f)(\kappa))$ is contained in $M$). It is shown in \cite{sy-glc}
that if $\delta$ is Woodin in $V=L[E]$ then in $L[E]$ there is a witness $T$ to 
Woodinness in this latter sense which can be lifted by a
non-set-generic real $R$, thereby preserving the Woodinness of
$\delta$. And indeed this can be done simultaneously for all Woodin cardinals in $L[E]$.

\bs

The proof of the latter result is much more involved than in the case
of nonstationary limits of strongs. In that simpler setting, one can
use the strength function $\alpha\mapsto (\sup$ of the strengths of
total extenders with critical point $\alpha)$ to thin out the
codings uniformly below each inaccessible cardinal. In the Woodin
cardinal setting one must instead use a uniform witness $\cal T$ to the Woodinness 
of each Woodin cardinal whose total extenders have non-Woodin critical
point, and then thin out the codings using functions which witness the
failure of these critical points to be Woodin. A major difference from
the easier setting is that for total extenders $F$ that are to be lifted 
and conditions $p$, it is no longer the case that $j_F(p)$ will be
trivial between the critical point and strength of $F$; instead one
must deal with this extra information at a cardinal $\alpha$ between
the critical point and strength of $F$ until reaching a condition
which ``recognises'' that each of the finitely-many total extenders in
$\cal T$ overlapping $\alpha$ has non-Woodin critical point; this is essential
for showing that this extra information stabilises to a set in $V$.

\bs

\noi
\emph{Future work}

\bs

The story is far from over regarding coding over core models. In terms
of versions of Corollary \ref{non-set-generic}, the current frontier
is the preservation of measurable Woodinness, which will 
need a technique beyond what is sketched above for plain Woodinness. 
Going back all the way to hypermeasurable cardinals, there remains the
difficult problem of condensation, which obstructs a satisfying
version of Theorem \ref{jensen-coding}. As mentioned, the special case
of coding a generic for a Prikry product is handled in
\cite{sy-mohammad}, but this is an extremely special case and it is
quite possible that there is a counterexample for the coding of more
general predicates while preserving hypermeasurability. And of course
it will be worthwhile to look at generalisations to the large cardinal
setting of the many applications of Jensen coding (and its
iterations), as found in \cite{sy-book, sy-david}. Finally, can
one do something with coding at the level of supercompact cardinals? Of course the core model
theory is not yet available there, but there has been
considerable progress in showing that many of the nice features of
$L[E]$ models can be forced consistently with the strongest of large
cardinal properties 
(see for example \cite{sy-peter}). Are there coding theorems to be
proved over such ``pseudo'' core models? A positive answer may have
very interesting consequences.

\end{document}